\documentclass[12pt]{amsart}
\usepackage{amsmath,amssymb}
\usepackage{amsthm}
\usepackage{thmtools}
\usepackage{mathtools}
\usepackage{dsfont}
\usepackage{color}
\usepackage{enumitem}
\usepackage{hyperref}
\usepackage[inner=1.0in,outer=1.0in,bottom=1.0in, top=1.0in]{geometry}
\usepackage{mleftright}
\usepackage{mathrsfs}
\usepackage{enumitem}
\usepackage{amssymb}

\usepackage{amsmath,amssymb}
\usepackage{tikz}
\usetikzlibrary{decorations.pathreplacing}
\usepackage{ytableau}

\usepackage{tikz}
\usetikzlibrary{decorations.pathreplacing}

\mleftright

\newcommand{\rk}{\operatorname{rank}}

\usepackage{nccmath}
\usepackage{cases}

\numberwithin{equation}{section}

\newtheorem{theorem}{Theorem}
\newtheorem*{example}{Example}

\newtheorem{proposition}[theorem]{Proposition}

\theoremstyle{remark}

\author[]{Colin Albert, Olivia Beckwith, Irfan Demetoglu, \\ Robert Dicks, John H. Smith, Jasmine Wang}
\title[]{
Integer partitions with large Dyson rank
}
\date{\today}

\begin{document}
\maketitle
\date{\today}

\begin{abstract}
The Dyson rank of an integer partition is the difference between its largest part and the number of parts it contains. Using Fine-Dyson symmetry,  we study counts of partitions whose rank lies in fixed residue classes and has large absolute value. We prove formulas relating these counts for partitions of different sizes. \end{abstract}

\section{Introduction}
A partition of an integer $n$ is a nonincreasing sequence of positive integers whose sum is $n$.  The partition function $p(n)$ is defined as the number of partitions of $n$ for $n\ge 1$, and $p(0)=1$. Throughout we let $P(n)$ denote the set of partitions of $n$, so that $p(n) = |P(n)|$ for integers $n \ge 1$, and $P(n) = \emptyset$ if $n$ is not a positive integer.

The arithmetic properties of the partition function have fascinated number theorists since Srinivasa Ramanujan \cite{Ram} discovered the Ramanujan congruences:
\begin{align*}
p(5k+4) &\equiv 0 \pmod{5}, \\
p(7k+5) &\equiv 0 \pmod{7}, \\
p(11k+6) &\equiv 0 \pmod {11}. 
\end{align*}

Ramanujan's proof of the Ramanujan congruences used properties of the generating function of $p(n)$. Freeman Dyson \cite{Dyson} conjectured that the first two Ramanujan congruences could be interpreted combinatorially using the \emph{rank statistic}, which he defined to be the difference between the largest element of the sequence and the number of elements of the sequence. He also conjectured that the third congruence could be explained combinatorially using a similar statistic called the \emph{crank}, which he did not claim to know how to define. His conjectures about the rank were proved by Atkin and Swinnerton-Dyer \cite{AtkinSwD}, and his conjectured crank function was discovered by Andrews and Garvan \cite{AndrewsGarvan}.

To state Dyson's conjectures about the rank more precisely, we let $\lambda = (\lambda_1, \lambda_2,\ldots,\lambda_k )$ denote a partition of $n$, with $\lambda_1 \ge \lambda_2 \ge\ldots\ge \lambda_k \ge 1$.  The rank of $\lambda$ is defined to be: 
\[\rk(\lambda):=\lambda_1-k\]
For example,
$$ \rk((6,3,3,2)) = 6-4 =2.$$

Dyson conjectured that the rank function divides the partitions of $5n+4$ (resp., $7n+5$) into $5$ (resp., $7$) different sets of equal size. More specifically, if we let $N(r,m,n)$ denote the number of partitions of $n$ whose rank is congruent to $r$ modulo $m$, then Dyson conjectured that for any $a,n \in \mathbb{Z}$,
$$
N(a, 5, 5n+4) = \frac{ p(5n+4)}{5}
$$
and similarly,
$$
N(a,7,7n+5) = \frac{p(7n+5)}{7}, 
$$
and since $N(a,5,5n+4)$ and $N(a,7,7n+5)$ are integers, these two claims imply the first and second Ramanujan congruences.

Dyson also conjectured several relations between counts of partitions with specific ranks, all of which were shown in \cite{AtkinSwD}. These relations are all of the form $N(a,p,pn+b) = N(c,p, pn+b)$ for $p \in \{5,7\}$, except for the following (see (19) of \cite{Dyson}): 
\begin{equation}\label{eq:dyson1}
N(0,7,7n+6) + N(1,7, 7n+6) = N(2,7,7n+6) + N(3,7,7n+6). 
\end{equation}
A few examples of relations similar to (\ref{eq:dyson1}) for ranks in residue classes modulo 8,9, and 12 are known by the work of Lewis \cite{Lewis} and Santa-Gadea \cite{SantaGadea}. We show that restricting to partitions such that the absolute value of the rank is larger than $\frac{n}{2}$, one obtains many identities resembling (\ref{eq:dyson1}). 

For any $n \in \mathbb{Z}$ we define $N^+(r,m,n)$ as follows:
$$
N^+(r,m,n) := \# \left\{\lambda \in P(n): \rk(\lambda) \equiv r \pmod{m}, |\rk(\lambda)| \geq \frac{n}{2}  \right\}  .
$$
Then we obtain the following family of identities:
\begin{theorem}\label{thm:resclasses}
Let $r_1,...r_4$, $a_1,...,a_4,m \in \mathbb{Z}$ be such that  $a_1,...,a_4,m$ are positive and we have:
\begin{align}
a_1-r_1 - m\Bigl\lceil\frac{a_1-2r_1}{2m}\Bigr\rceil &=a_3-r_3-m\Bigl\lceil\frac{a_3-2r_3}{2m}\Bigr\rceil, \label{eq:condition1} \\
a_1+r_1- m\Bigl\lceil\frac{a_1+2r_1}{2m}\Bigr\rceil & =a_4-r_4-m\Bigl\lceil\frac{a_4-2r_4}{2m}\Bigr\rceil,\label{eq:condition2} \\
a_2-r_2-m\Bigl\lceil\frac{a_2-2r_2}{2m}\Bigr\rceil &=a_4+r_4-m\Bigl\lceil\frac{a_4+2r_4}{2m}\Bigr\rceil,\label{eq:condition3}  \\
a_2+r_2-m\Bigl\lceil\frac{a_2+2r_2}{2m}\Bigr\rceil&=a_3+r_3-m\Bigl\lceil\frac{a_3+2r_3}{2m}\Bigr\rceil \label{eq:condition4} .
\end{align}
Then for $n \ge 0$, we have:
\begin{equation*}
N^+(r_1,m,2mn+a_1) + N^+(r_2,m,2mn+a_2) = N^+(r_3,m,2mn+a_3)+N^+(r_4,m,2mn+a_4).
\end{equation*}
\end{theorem}
Note that if the $a_i$ are all equal, one has $r_1 \equiv r_3 \equiv -r_4 \equiv -r_2 \equiv -r_3 \pmod{m}$, and the statement is not interesting. However, when the $a_i$ are distinct, one obtains relations between counts of partitions of different numbers with rank that have large absolute values and lie in given residue classes. Comparing counts of partitions of different numbers with prescribed rank properties is significantly different than (\ref{eq:dyson1}).

\begin{example} Theorem \ref{thm:resclasses} implies that for $n \ge 0$,
\begin{align*}
N^+(2,19,38n+20) + N^+(5,19,38n+22) &= N^+(14,19,38n+32)+N^+(7,19,38n+29) \\
N^+(1,21,42n+43) + N^+ (19,21,42n+50) &= N^+(3,21,42n+45) + N^+ (4,21,42n+27) \\
N^+(1,23,46n+47) + N^+(-2,23,46n+54) &= N^+(3,23,46n+49) + N^+(4,23,46n+29).
\end{align*}
\end{example}
It is not difficult in practice to find $m, r_i, a_i$ satisfying (\ref{eq:condition1})-(\ref{eq:condition4}). Given a choice of $a_1, \cdots, a_4$ and $r_1, \cdots, r_4$ for which (\ref{eq:condition1})-(\ref{eq:condition4}) hold modulo $m$, it is usually (although not always) possible to shift the $a_i$ modulo $2m$, thereby shifting each side of (\ref{eq:condition1})-(\ref{eq:condition4}) by a multiple of $m$, to produce a solution of (\ref{eq:condition1})-(\ref{eq:condition4}). More examples can be found in the Mathematica \cite{math} code on the second author's homepage. 

Theorem \ref{thm:resclasses} is proved using Fine-Dyson symmetry (see \cite{Fine}, \cite{Fine2},  \cite{Dyson2}, and the survey \cite{Pak}), which relates the partitions of $n$ with rank at most $r$ to partitions of $n+r-1$ with rank at most $2-r$ using a bijective argument on Young diagrams. As a consequence, if one lets $N(m,n)$ denote the number of partitions of $n$ with rank $m$, then one can show that the inequality \begin{equation}\label{eq:Prop2.1} N(m,n) \le p(n-1-|m|) - p(n - 2 - |m|), \end{equation} holds, with equality if $|m| \ge \frac{n}{2} - 2$. 

In Section 2 we give an alternative proof of this fact. Our proof involves the set of partitions for which the number of parts is a part. We let $g(k)$ denote the number of partitions of $k$ satisfying this property.
\begin{proposition}\label{thm:symmetry}
For all $m,n \in \mathbb{Z}$, we have
$$ N(m,n) \le g( n - |m|) = p(n-1-|m|) - p(n - 2 - |m|),$$
with equality if $|m| \ge \frac{n}{2} - 2$.
\end{proposition}
Here we use the convention that $p(x) = 0$ if $x$ is not a non-negative integer. 


Considering the partitions of $n$ with large rank was partly inspired by asymptotic formulas for the number of partitions of $n$ with rank or crank equal to $m$, where $|m| \le \frac{\sqrt{n} \log{n}}{\pi \sqrt{6}} $. These formulas were conjectured by Dyson \cite{Dyson3} and proved using the circle method by Bringmann and Dousse \cite{BringmannDousse} for cranks and by Dousse and Mertens \cite{DousseMertens} for ranks.  In contrast to the partitions with small rank, the partitions considered in this paper are easy to count: Proposition \ref{thm:symmetry} implies that the number of partitions of $n$ with rank at least $\frac{n}{2}$ is $p(\lfloor \frac{n}{2} \rfloor-1)$. By the Hardy-Ramanujan asymptotic \cite{HardyRam},  we have
$$p(n)\sim\frac{1}{4n\sqrt{3}}e^{\pi\sqrt{\frac{2n}{3}}}$$
Thus, the portion of the partitions that have rank at least $n/2$ is exponentially decreasing with $n$: 
$$
\frac{ \# \{ \lambda \in P(n): \rk(\lambda) \ge \frac{n}{2} \}}{p(n)} = \frac{p(\lfloor \frac{n}{2} \rfloor -1) }{p(n)} \sim 2e^{\pi(1-\sqrt{2})\sqrt{\frac{n}{3}}}.$$
For example, there are 37338 partitions of 40 (shown in the histogram), but only 490 have rank greater than or equal to 20, a portion of about 1.31\%. 
\begin{figure}
\caption{Ranks of partitions $\lambda \in P(40)$}
\includegraphics[width=5cm]{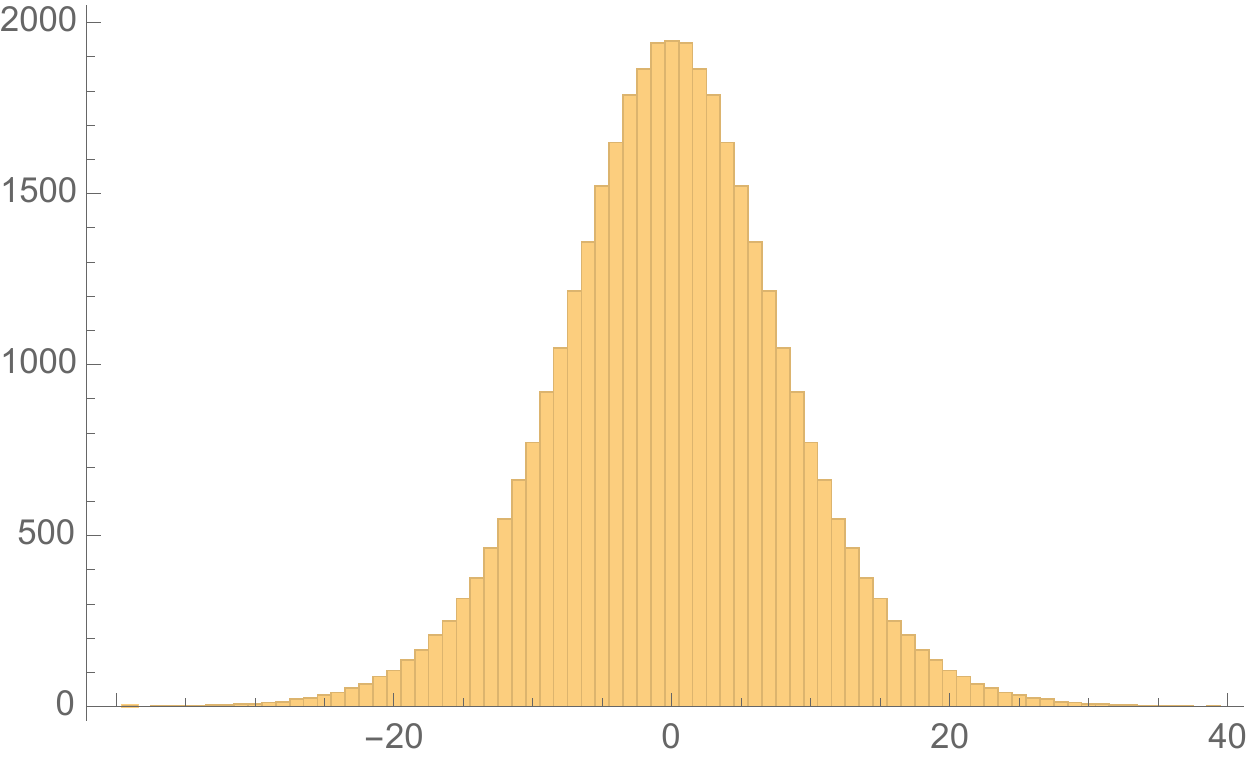}
\end{figure}

In Section 2 we review some standard facts about partitions and prove Proposition \ref{thm:symmetry}. The proof of Theorem \ref{thm:resclasses} is given in Section 3. Section 4 mentions a few questions for future work.

\section*{Acknowledgements}
The authors are grateful to the anonymous referee for a thorough report offering several valuable suggestions that have improved the manuscript. 
The authors also thank the Illinois Geometry Lab, which is supported by the Department of Mathematics at the University of Illinois at Urbana-Champaign. This material is based upon work supported by the National Science Foundation under Grant No. DMS-1449269. Any opinions, findings, and conclusions or recommendations expressed in this material are those of the authors and do not necessarily reflect the views of the National Science Foundation.

\section{Fine-Dyson Symmetry}
Throughout, we let $\lambda = (\lambda_1,\ldots,\lambda_b)$ and $\mu = (\mu_1,\ldots,\mu_c)$ denote partitions, with the $\lambda_i$ and $\mu_i$ positive integers satisfying $\lambda_1 \ge \lambda_2 \ge\ldots\ge\lambda_b$ and $\mu_1 \ge \mu_2 \ge \ldots\ge \mu_c$. The $\lambda_i$ are called the \emph{parts} of $\lambda$. We let $|\lambda| = \lambda_1 + \lambda_2 +\ldots+ \lambda_b$, and we let $\lambda \cup \mu$ be the partition of $|\lambda| + |\mu|$ consisting of the union of the parts of both $\lambda$ and $\mu$, e.g.
$$
(4,2,2,1) \cup (7,4,2) = (7,4,4,2,2,2,1).
$$

Partitions are naturally visualized using Young diagrams. The Young diagram of a partition $\lambda$ is an array of squares such that the top row has $\lambda_1$ squares, the second row has $\lambda_2$ squares, and so on. For example, the Young diagram for $\lambda = (3,2,2,1)$ is

\begin{center}
\begin{ytableau}
\none      &  &  &   \\
\none   & &   \\
\none & &  \\
\none & 
\end{ytableau}
\end{center}
Reflecting the Young diagram across the $y=-x$ line, where the top left corner is the origin $(0,0)$, produces a Young diagram for another partition, called the \emph{conjugate} of $\lambda$. For example, the conjugate of $\lambda = (3,2,2,1)$ is $(4,3,1)$: 

\begin{center}
\begin{ytableau}
\none      &  &  & &  \\
\none   & &  & \\
\none &  \\
\end{ytableau}
\end{center}

The rank of the conjugate of a partition $\lambda$ is $-\rk(\lambda)$, and it follows that $N(m,n) = N(-m,n)$.


\subsection{Proof of Proposition \ref{thm:symmetry}}
The claim (\ref{eq:Prop2.1}) follows from work of Fine and Dyson  (\cite{Fine}, \cite{Fine2},  \cite{Dyson2}), and we'll briefly summarize their techniques before presenting our alternative proof. 

Let $h(n,r)$ denote the number of partitions of $n$ with rank at most $r$, thus, 
$$
h(n,r) = \sum_{m=-(n-1)}^{r} N(m,n).
$$ 
Fine and Dyson define a bijection between the set of partitions of $n$ of rank at most $r$ and the set of partitions of $n + r -1$ of rank at least $r-2$, in which one adds $r-1$ squares to the leftmost column of a Young diagram of a partition of $n$ and then moves it to the top row. Conjugating the partitions of $n + r -1$ with rank at least $r-2$, one obtains
\begin{equation}\label{eq:FineDyson}
h(n,r) = h(n+r-1,2-r).
\end{equation}

Then (\ref{eq:Prop2.1}) can be derived from (\ref{eq:FineDyson}) and $N(m,n) = h(n,m) - h(n,m-1)$. 

\begin{proof}
 First we will show that $N(m,n) \le g(n-|m|)$, with equality if $m \ge \frac{n}{2} - 2$. 

From the rank symmetry $N(m,n) = N(-m,n)$, we may assume that $m$ is non-negative without loss of generality.

We will define a map $f$ from the partitions of $n$ with rank $m$ to the partitions of $n-m$ such that the number of parts is a part. Take a partition $\lambda$ of $n$ with rank $m$. Then it is of the form:
\begin{equation*}
  \lambda =   (m + c) \cup \mu
\end{equation*}
where $m+c$ is the largest part of $\lambda$ and $c$ is the number of parts of $\lambda$, with $\mu$ a partition of $n-m-c$. Then we define $f(\lambda)$ to be $f(\lambda):= (c) \cup \mu.$ We note that $f(\lambda)$ is a partition of $n-m$ with $c$ parts, and with $c$ present as a part. 

Suppose $\lambda = (m+c, \lambda_2,\ldots,\lambda_c)$ and $\lambda' = (m+c', \lambda_2',\ldots, \lambda_{c'}')$ are such that
$$
f(\lambda) = \mu  = f(\lambda'),
$$
where $\mu$ is a partition of $n-m$. Then $c$ and $c'$ both equal the number of parts of $\mu$, and hence $c'=c$. It then follows from $f(\lambda) = f(\lambda')$ that $\lambda_i = \lambda_i'$ for each $1 \le i \le c$, so that $\lambda = \lambda'$. Thus the map $f$ is injective, which proves the first part of the claim, that $N(m,n) \le g(n-|m|)$. 

To show the second half of the claim, we assume $m \ge \frac{n}{2}-2$ . We will show that $f$ is also surjective: first, if $m=n-1$, then the only partition of $n$ with rank $m$ is the partition $(n)$, and this corresponds to the only partition of $n-m=1$, for which the number of parts is equal to one, and is the size of the only part. Now suppose $m \le n-2$, and let $\mu$ be a partition of $n-m$ whose number of parts is a part. Then $\mu$ is of the form $(c) \cup \mu'$, where $\mu'$ is a partition of $n-m-c$ with $c-1$ parts. Note that $c \ge 2$, since $n-m \ge 2$. We observe that $\delta := (c+m) \cup \mu'$ is a partition of $n$ with $c$ parts. The largest part of $\delta$ must be $c+m$, since $$c+m \ge m+2 \ge \frac{n}{2} \ge n-m-c = |\mu'|.$$

Finally, we have \begin{equation}\label{eq:bijection2}g(k)=p(k-1)-p(k-2)\end{equation} for all $k \ge 1$  (see the sequence A002865 \cite{oeis}). For completeness, we provide a brief proof. Clearly, $p(k-1) - p(k-2)$ is the number of partitions of $k-1$ for which one is not a part. Given such a partition $\lambda = (\lambda_1,\ldots,\lambda_b)$, we can construct a partition $\gamma = (\lambda_1 - 1, \lambda_2 - 1,\ldots,\lambda_b -1) \cup (b+1)$, which is partition of $k$ with $b+1$ parts. On the other hand, if $\mu = (\mu_1,\ldots, \mu_b)$ is a partition of $k$ such that $b$ is a part, then $(\mu_1 + 1,\ldots,\mu_b + 1)\backslash (b + 1)$ is a partition of $k-1$ which contains no ones as parts. We have given a bijection between the partitions of $k-1$ containing no parts equal to one and the partitions of $k$ for which the number of parts is a part, which proves (\ref{eq:bijection2}).

Thus we have shown that 
$$N(m,n) \le g(n-|m|) = p(n-|m|-1) - p(n-|m| - 2),$$
 with equality if $|m| \ge \frac{n}{2} -2$.
\end{proof}
\section{Proof of Theorem \ref{thm:resclasses}} 
    Assume that $m, a_1, \cdots , a_4, r_1, \cdots , r_4$ satisfy the assumptions of Theorem \ref{thm:resclasses}. By definition, 
$$
N^+(r_j,m,2mn+a_j) = \sum_{|r_j+mk| \geq \frac{2mn+a_j}{2} } N(r_j+mk,2mn+a_j).
$$
Applying Proposition \ref{thm:symmetry}, we obtain
\begin{align*}
N^+(r_j,m,2mn+a_j) & = \sum_{r_j+mk \geq \frac{2mn+a_j}{2} }N(r_j+mk,2mn+a_j) \\
&\quad \quad + \sum_{r_j+mk \leq \frac{-(2mn+a_j)}{2} }N(r_j+mk,2mn+a_j) \\
&=\sum_{k \geq n+ \lceil \frac{a_j-2r_j}{2m} \rceil  }g(a_j-r_j+m(2n-k)) \\
&\quad \quad+ \sum_{k \geq n +  \lceil \frac{a_j+2r_j}{2m}\rceil  }g(a_j+r_j+m(2n-k)) \\
&= S_j^- + S_j^+,
\end{align*}
where
\begin{equation}\label{eq:Sminus}
S_j^- := \sum_{k \geq n + \lceil \frac{a_j-2r_j}{2m} \rceil  }g(a_j-r_j+m(2n-k)) 
\end{equation}
and
\begin{equation}\label{eq:Splus}
S_j^+ :=  \sum_{k \geq n +  \lceil \frac{a_j+2r_j}{2m}\rceil  }g(a_j+r_j+m(2n-k)).
\end{equation}

With this notation, our goal is to show
$$
    S_1^{+} + S_1^{-} + S_2^{+} + S_2^{-} = S_3^{+} + S_3^{-} + S_4^{+} + S_4^{-}.
$$
We will prove this by showing the following equalities:
\begin{align*}
S_1^{-} &= S_3^{-}, \\
S_1^{+} &= S_4^{-}, \\
S_2^{-} &= S_4^{+}, \\
 S_2^{+} &= S_3^{+}. \\
\end{align*}

Let $t_{ij}:=\frac{a_i-r_i-(a_j-r_j)}{m}$, then we have:
$$
S_i^- =
\sum_{k \geq n +\lceil \frac{a_i-2r_i}{2m} \rceil  } g(a_j-r_j+m(2n-k+t_{ij})).
$$
Next we make the substitution $v=k-t_{ij}$ and obtain:
$$
S_i^- = \sum_{v \geq n - t_{ij} +\lceil \frac{a_i-2r_i}{2m} \rceil } g(a_j-r_j+m(2n-v)).
$$
This is equal to $S_j^-$ if 
$$
n - t_{ij} + \left \lceil \frac{a_i-2r_i}{2m} \right \rceil = n +\left \lceil \frac{a_j-2r_j}{2m} \right \rceil.
$$
Using (\ref{eq:condition1}), it follows that $S_1^- = S_3^-.$ 

Let $s_{ij}:=\frac{a_i-r_i-(a_j+r_j)}{m}$, then we have:
$$
S_i^-=
\sum_{k \geq n +\lceil \frac{a_i-2r_i}{2m} \rceil } g(a_j+r_j+m(2n-k+s_{ij})).
$$
Next we make the substitution $v=k-s_{ij}$ and obtain:
$$
S_i^- = \sum_{v \geq n +\lceil \frac{a_i-2r_i}{2m} \rceil - s_{ij} } g(a_j+r_j+m(2n-v)).
$$
This is equal to $S_j^+$ if $\lceil \frac{a_i-2r_i}{2m} \rceil - s_{ij} =  n +  \lceil \frac{a_j+2r_j}{2m}\rceil$. Therefore $S_1^{+} = S_4^{-}$ and  $S_2^{-} = S_4^{+}$ follow from (\ref{eq:condition2})-(\ref{eq:condition3}).

The last equation, $S_2^+ = S_3^+$ follows almost identically: setting $w = \frac{a_2 + r_2 - (a_3 + r_3)}{m}$, we have
$$
S_2^+ = \sum_{k \geq n +  \lceil \frac{a_2+2r_2}{2m}\rceil }g(a_3+r_3+m(2n-k + w) ) =  \sum_{v \geq n - w +  \lceil \frac{a_2+2r_2}{2m}\rceil  }g(a_3+r_3+m(2n-v) ). 
$$
This is equal to $S_3^+$ if 
$$
n - w +  \left\lceil \frac{a_2+2r_2}{2m} \right\rceil  = n +  \left\lceil \frac{a_3+2r_3}{2m}\right\rceil 
$$
which follows from (\ref{eq:condition4}).

\section{Future work}
The combinatorial nature of the identities in Theorem \ref{thm:resclasses} raises the question of whether a combinatorial proof could be found that directly relates the partitions of the $a_i + 2 m n$, as opposed to relying on the connection Fine-Dyson symmetry provides between the partitions of $n$ with rank $m$ and the partitions of $n-|m|$.  

In a separate direction, it would be interesting to find analogues of Theorem \ref{thm:resclasses} for other partition statistics, such as the crank. 

Finally, it is natural to wonder what more can be said about the distribution of ranks.  As mentioned in the introduction, an immediate consequence of (\ref{eq:Prop2.1}) is that the number of partitions with rank at least $\frac{n}{2}$ is $p(\lfloor \frac{n}{2} \rfloor-1)$. More generally, (\ref{eq:Prop2.1}) implies a formula for the number of partitions with rank in intervals contained in $[n/2,n]$: if $b \ge a \ge \frac{n}{2} - 2$, then
$$ \left| \{\lambda \in P(n) :  \rk(\lambda) \in [a,b] \} \right| = p(n-\lceil a \rceil-1)-p(n-\lfloor b \rfloor-2).$$
Using completely different techniques (the circle method), a few aforementioned works (\cite{Dyson3}, \cite{BringmannDousse}, \cite{DousseMertens}) provide asymptotic formulas for $N(m,n)$ with $|m| \le \frac{\sqrt{n} \log{n}}{\pi \sqrt{6}} $. Little is known, to our knowledge, about the number of partitions with rank in ranges between $\frac{\sqrt{n} \log{n}}{\pi \sqrt{6}} $ and $n/2$, and it would be very interesting to obtain formulas for these counts.



  \bibliographystyle{alpha}
\bibliography{final.bib} 

\begin{thebibliography}{ASD54}

\bibitem[AG88]{AndrewsGarvan}
G.~E. Andrews and F.~G. Garvan.
\newblock Dyson's crank of a partition.
\newblock {\em Bull. Amer. Math. Soc. (N.S.)}, 18(2):167--171, 1988.

\bibitem[ASD54]{AtkinSwD}
A.O.L. Atkin and P.~Swinnerton-Dyer.
\newblock Some properties of partitions.
\newblock {\em Proc. London Math. Soc. (3)}, 4(2), 1954.

\bibitem[BD16]{BringmannDousse}
K.~Bringmann and J.~Dousse.
\newblock On {D}yson's crank conjecture and the uniform asymptotic behavior of
  certain inverse theta functions.
\newblock {\em Trans. Amer. Math. Soc.}, 368(5):3141--3155, 2016.

\bibitem[DM15]{DousseMertens}
J.~Dousse and M.~H. Mertens.
\newblock Asymptotic formulae for partition ranks.
\newblock {\em Acta Arith.}, 168(1):83--100, 2015.

\bibitem[Dys44]{Dyson}
F.J. Dyson.
\newblock Some guesses in the theory of partitions.
\newblock {\em Eureka (Cambridge)}, 8, 1944.

\bibitem[Dys69]{Dyson2}
F.~J. Dyson.
\newblock A new symmetry of partitions.
\newblock {\em J. Combinatorial Theory}, 7:56--61, 1969.

\bibitem[Dys89]{Dyson3}
F.~J. Dyson.
\newblock Mappings and symmetries of partitions.
\newblock {\em J. Combin. Theory Ser. A}, 51(2):169--180, 1989.

\bibitem[Fin48]{Fine}
N.~J. Fine.
\newblock Some new results on partitions.
\newblock {\em Proc. Nat. Acad. Sci. U.S.A.}, 34:616--618, 1948.

\bibitem[Fin88]{Fine2}
N.~J. Fine.
\newblock {\em Basic hypergeometric series and applications}, volume~27 of {\em
  Mathematical Surveys and Monographs}.
\newblock American Mathematical Society, Providence, RI, 1988.
\newblock With a foreword by George E. Andrews.

\bibitem[HR00]{HardyRam}
G.~H. Hardy and S.~Ramanujan.
\newblock Asymptotic formul\ae in combinatory analysis [{P}roc. {L}ondon
  {M}ath. {S}oc. (2) {\bf 16} (1917), {R}ecords for 1 {M}arch 1917].
\newblock In {\em Collected papers of {S}rinivasa {R}amanujan}, page 244. AMS
  Chelsea Publ., Providence, RI, 2000.

\bibitem[Lew92]{Lewis}
R.~Lewis.
\newblock On some relations between the rank and the crank.
\newblock {\em J. Combin. Theory Ser. A}, 59(1):104--110, 1992.

\bibitem[oei]{oeis}
Oeis foundation inc. (2021), the on-line encyclopedia of integer sequences,
  https://oeis.org/.

\bibitem[Pak06]{Pak}
I.~Pak.
\newblock Partition bijections, a survey.
\newblock {\em Ramanujan J.}, 12(1):5--75, 2006.

\bibitem[Ram21]{Ram}
S.~Ramanujan.
\newblock Congruence properties of partitions.
\newblock {\em Math. Z.}, 9(1-2):147--153, 1921.

\bibitem[SG90]{SantaGadea}
N.~A. Santa~Gadea.
\newblock {\em On the rank and the crank moduli 8, 9 and 12}.
\newblock ProQuest LLC, Ann Arbor, MI, 1990.
\newblock Thesis (Ph.D.)--The Pennsylvania State University.

\bibitem[{Wol}]{math}
{Wolfram Research, Inc.}
\newblock Mathematica 8.0.

\end{thebibliography}
\end{document}